\documentclass[12pt]{amsart}

\usepackage{amsmath}
\usepackage{amsxtra}
\usepackage{amscd}
\usepackage{amsthm}
\usepackage{amsfonts}
\usepackage{amssymb}
\usepackage{eucal}
\usepackage[matrix,arrow,curve]{xy}
\usepackage[dvips]{graphicx}
\usepackage[dvips]{graphics}
\usepackage[T2A]{fontenc}
\usepackage[cp866]{inputenc}
\usepackage{mathtools}

\newcommand{\End}{\mathop{\sf End}\nolimits}
\newcommand{\Sym}{\mathop{\textrm{Sym}}\nolimits}

\newcommand{\Ker}{\mathop{\sf Ker}\nolimits}
\newcommand{\im}{\mathop{\sf Im}\nolimits}

\newcommand{\Sec}{\mathop{\textrm{Sec}}\nolimits}

\def \id {{\rm id}}

\theoremstyle{plain}
\newtheorem{Thm}[subsection]{Theorem}
\newtheorem{Cor}[subsection]{Corollary}
\newtheorem{Lem}[subsection]{Lemma}
\newtheorem{Prop}[subsection]{Proposition}
\newtheorem{Conj}[subsection]{Conjecture}
\newtheorem{Ex}[subsection]{Example}

\theoremstyle{definition}
\newtheorem{Def}[subsection]{Definition}

\theoremstyle{remark}

\newtheorem{Rem}[subsection]{Remark}

\newtheorem{conj}{Conjecture}

\numberwithin{equation}{section}

\newif\ifShowLabels
\ShowLabelstrue
\newdimen\theight
\def\TeXref#1{%
    \leavevmode\vadjust{\setbox0=\hbox{{\tt
        \quad\quad  {\small \rm #1}}}%
    \theight=\ht0
    \advance\theight by \lineskip
    \kern -\theight \vbox to
    \theight{\rightline{\rlap{\box0}}%
    \vss}%
    }}%

\ShowLabelsfalse

\renewcommand{\sec}[2]{\section{#2}\label{S:#1}%
    \ifShowLabels \TeXref{{S:#1}} \fi}
\newcommand{\ssec}[2]{\subsection{#2}\label{SS:#1}%
    \ifShowLabels \TeXref{{SS:#1}} \fi}

\newcommand{\refss}[1]{Section ~\ref{SS:#1}}

\newcommand{\reft}[1]{Theorem ~\ref{T:#1}}
\newcommand{\refl}[1]{Lemma ~\ref{L:#1}}

\newenvironment{thm}[1]%
    { \begin{Thm} \label{T:#1}  \ifShowLabels \TeXref{T:#1} \fi }%
    { \end{Thm} }

\renewcommand{\th}[1]{\begin{thm}{#1} \sl }
\renewcommand{\eth}{\end{thm} }

\newenvironment{lemma}[1]%
    { \begin{Lem} \label{L:#1}  \ifShowLabels \TeXref{L:#1} \fi }%
    { \end{Lem} }
\newcommand{\lem}[1]{\begin{lemma}{#1} \sl}
\newcommand{\elem}{\end{lemma}}

\newenvironment{propos}[1]%
    { \begin{Prop} \label{P:#1}  \ifShowLabels \TeXref{P:#1} \fi }%
    { \end{Prop} }
\newcommand{\prop}[1]{\begin{propos}{#1}\sl }
\newcommand{\eprop}{\end{propos}}

\newenvironment{corol}[1]%
    { \begin{Cor} \label{C:#1}  \ifShowLabels \TeXref{C:#1} \fi }%
    { \end{Cor} }
\newcommand{\cor}[1]{\begin{corol}{#1} \sl }
\newcommand{\ecor}{\end{corol}}

\newenvironment{defeni}[1]%
    { \begin{Def} \label{D:#1}  \ifShowLabels \TeXref{D:#1} \fi }%
    { \end{Def} }
\newcommand{\defe}[1]{\begin{defeni}{#1} \sl }
\newcommand{\edefe}{\end{defeni}}

\newenvironment{remark}[1]%
    { \begin{Rem} \label{R:#1}  \ifShowLabels \TeXref{R:#1} \fi }%
    { \end{Rem} }
\newcommand{\rem}[1]{\begin{remark}{#1}}
\newcommand{\erem}{\end{remark}}

\newenvironment{conjec}[1]%
    { \begin{Conj} \label{Co:#1}  \ifShowLabels \TeXref{Co:#1} \fi }%
    { \end{Conj} }
\renewcommand{\conj}[1]{\begin{conjec}{#1} \sl }
\newcommand{\econj}{\end{conjec}}

\newenvironment{example}[1]%
    { \begin{Ex} \label{Exx:#1}  \ifShowLabels \TeXref{Exx:#1} \fi }%
    { \end{Ex} }
\newcommand{\ex}[1]{\begin{example}{#1} \sl }
\newcommand{\eex}{\end{example}}

\newcommand{\eq}[1]%
    { \ifShowLabels \TeXref{E:#1} \fi
       \begin{equation} \label{E:#1} }
\newcommand{\eeq}{ \end{equation} }

\newcommand{\prf}{ \begin{proof} }
\newcommand{\epr}{ \end{proof} }

\newcommand\alp{\alpha}     
\newcommand\bet{\beta}
     \newcommand\Gam{\Gamma}
     \newcommand\Del{\Delta}
\newcommand\eps{\varepsilon}

\newcommand\iot{\iota}

\newcommand\lam{\lambda}        \newcommand\Lam{\Lambda}
     \newcommand\Sig{\Sigma}

\newcommand\ome{\omega}     

\newcommand\calC{{\mathcal{C}}}

\newcommand\calW{{\mathcal{W}}}

        \newcommand\bfM{{\mathbf M}}

\newcommand\LL{\mathbb{L}}

 \newcommand\grg{{\mathfrak{g}}}

\newcommand\sdp{\times \hskip -0.3em {\raise 0.3ex
\hbox{$\scriptscriptstyle |$}}} 

\newcommand\Coker{\operatorname{Coker}}

\newcommand\Ext{\operatorname{Ext}}

\newcommand\Gr{\operatorname{Gr}}

\newcommand\Int{\operatorname{Int}}

\newcommand\proj{\operatorname{proj}}

\newcommand\rk{\operatorname{rk}}

\newcommand\Perv{\operatorname{Perv}}

\newcommand\SL{{\rm SL}}

\newcommand\supp{\operatorname{supp}}

\newcommand\uc{{\underline{c}}}

\newcommand\tilc{{\widetilde{c}}}
\newcommand\tilC{{\widetilde{C}}}

\newcommand\tilH{{\widetilde{H}}}

\newcommand\tilS{{\widetilde{S}}}

\newcommand\tilZ{{\widetilde{Z}}}

\newcommand\tilchi{{\widetilde{\chi}}}

\newcommand\tilpi{{\widetilde{\pi}}}

\newcommand\x{\times}
\newcommand\ten{\otimes}

\newcommand\nc{\newcommand}

\newcommand{\IC}{{\operatorname{IC}}}

\newcommand{\disc}{\operatorname{disc}}

\newcommand{\iso}{{\stackrel{\sim}{\longrightarrow}}}

\nc\aff{\operatorname{aff}}
\nc\oGr{\overline{\Gr}}
\nc\Bun{\operatorname{Bun}}
\nc\hgrg{\widehat{\grg}}
\renewcommand\Int{\operatorname{Int}}
\nc\bInt{\overline{\Int}}
\nc\hatLam{\widehat{\Lam}}
\nc\bmu{\overline{\mu}}
\nc\bnu{\overline{\nu}}
\nc\blambda{\overline{\lam}}
\renewcommand\SL{\operatorname{SL}}
\nc\ocalW{\overline{\calW}}
\nc\pos{\operatorname{pos}}
\nc\IH{\operatorname{IH}}
\nc\Rep{\operatorname{Rep}}
\nc\Gal{\operatorname{Gal}}
\nc{\tilGr}{\widetilde{\Gr}}

\nc\Pic{\operatorname{Pic}}
\nc\Prym{\operatorname{Prym}}
\nc\pa{\partial}
\nc\Na{\nabla}

\nc{\HC}{{\mathcal{HC}}}
\nc{\on}{\operatorname}
\nc{\BA}{{\mathbb{A}}}
\nc{\BC}{{\mathbb{C}}}
\nc{\BG}{{\mathbb{G}}}
\nc{\BM}{{\mathbb{M}}}
\nc{\BN}{{\mathbb{N}}}
\nc{\BQ}{{\mathbb{Q}}}
\nc{\BP}{{\mathbb{P}}}
\nc{\BR}{{\mathbb{R}}}
\nc{\BZ}{{\mathbb{Z}}}
\nc{\BS}{{\mathbb{S}}}

\nc{\CA}{{\mathcal{A}}}
\nc{\CB}{{\mathcal{B}}}
\nc{\CalC}{{\mathcal C}}
\nc{\CalD}{{\mathcal D}}
\nc{\CE}{{\mathcal{E}}}
\nc{\CF}{{\mathcal{F}}}
\nc{\CG}{{\mathcal{G}}}
\nc{\CH}{{\mathcal{H}}}
\nc{\CK}{{\mathcal{K}}}
\nc{\CL}{{\mathcal{L}}}
\nc{\CM}{{\mathcal{M}}}
\nc{\CMM}{{\mathcal{M}^{\operatorname{gen}}_\hbar(-\rho)}}
\nc{\CN}{{\mathcal{N}}}
\nc{\CO}{{\mathcal{O}}}
\nc{\CP}{{\mathcal{P}}}
\nc{\CQ}{{\mathcal{Q}}}
\nc{\CR}{{\mathcal{R}}}
\nc{\CS}{{\mathcal{S}}}
\nc{\CT}{{\mathcal{T}}}
\nc{\CU}{{\mathcal{U}}}
\nc{\CV}{{\mathcal{V}}}
\nc{\CW}{{\mathcal{W}}}
\nc{\CX}{{\mathcal{X}}}
\nc{\CY}{{\mathcal{Y}}}
\nc{\CZ}{{\mathcal{Z}}}

\nc{\gen}{{\operatorname{gen}}}
\nc{\cM}{{\check{\mathcal M}}{}}
\nc{\csM}{{\check{\mathcal A}}{}}
\nc{\obM}{{\overset{\circ}{\mathbf M}}{}}
\nc{\oCA}{{\overset{\circ}{\mathcal A}}{}}
\nc{\obA}{{\overset{\circ}{\mathbf A}}{}}
\nc{\ooM}{{\overset{\circ}{M}}{}}
\nc{\osM}{{\overset{\circ}{\mathsf M}}{}}
\nc{\vM}{{\overset{\bullet}{\mathcal M}}{}}
\nc{\nM}{{\underset{\bullet}{\mathcal M}}{}}
\nc{\obD}{{\overset{\circ}{\mathbf D}}{}}
\nc{\cp}{{\overset{\circ}{\mathbf p}}{}}
\nc{\ofZ}{{\overset{\circ}{\mathfrak Z}}{}}

\nc{\fa}{{\mathfrak{a}}}
\nc{\fb}{{\mathfrak{b}}}
\nc{\fg}{{\mathfrak{g}}}
\nc{\fgl}{{\mathfrak{gl}}}
\nc{\fh}{{\mathfrak{h}}}
\nc{\fj}{{\mathfrak{j}}}
\nc{\fm}{{\mathfrak{m}}}
\nc{\fn}{{\mathfrak{n}}}
\nc{\fu}{{\mathfrak{u}}}
\nc{\fp}{{\mathfrak{p}}}
\nc{\frr}{{\mathfrak{r}}}
\nc{\fs}{{\mathfrak{s}}}
\nc{\ft}{{\mathfrak{t}}}
\nc{\fT}{{\mathfrak{T}}}
\nc{\ofT}{{\overline{\mathfrak T}}}
\nc{\ofS}{{\overline{\mathfrak S}}}
\nc{\fsl}{{\mathfrak{sl}}}
\nc{\hsl}{{\widehat{\mathfrak{sl}}}}
\nc{\hgl}{{\widehat{\mathfrak{gl}}}}
\nc{\hg}{{\widehat{\mathfrak{g}}}}
\nc{\chg}{{\widehat{\mathfrak{g}}}{}^\vee}
\nc{\hn}{{\widehat{\mathfrak{n}}}}
\nc{\chn}{{\widehat{\mathfrak{n}}}{}^\vee}

\nc{\fA}{{\mathfrak{A}}}
\nc{\fB}{{\mathfrak{B}}}
\nc{\fD}{{\mathfrak{D}}}
\nc{\fE}{{\mathfrak{E}}}
\nc{\fF}{{\mathfrak{F}}}
\nc{\fG}{{\mathfrak{G}}}
\nc{\fI}{{\mathfrak{I}}}
\nc{\fJ}{{\mathfrak{J}}}
\nc{\fK}{{\mathfrak{K}}}
\nc{\fL}{{\mathfrak{L}}}
\nc{\fM}{{\mathfrak{M}}}
\nc{\fN}{{\mathfrak{N}}}
\nc{\frP}{{\mathfrak{P}}}
\nc{\fS}{{\mathfrak S}}
\nc{\fU}{{\mathfrak{U}}}
\nc{\fZ}{{\mathfrak{Z}}}

\nc{\bb}{{\mathbf{b}}}
\nc{\bc}{{\mathbf{c}}}
\nc{\be}{{\mathbf{e}}}
\nc{\bj}{{\mathbf{j}}}
\nc{\bn}{{\mathbf{n}}}
\nc{\bp}{{\mathbf{p}}}
\nc{\bq}{{\mathbf{q}}}
\nc{\bv}{{\mathbf{v}}}
\nc{\bx}{{\mathbf{x}}}
\nc{\by}{{\mathbf{y}}}
\nc{\bw}{{\mathbf{w}}}
\nc{\bA}{{\mathbf{A}}}
\nc{\bB}{{\mathbf{B}}}
\nc{\bC}{{\mathbf{C}}}
\nc{\bK}{{\mathbf{K}}}
\nc{\bL}{{\mathbf{L}}}
\nc{\bD}{{\mathbf{D}}}
\nc{\bH}{{\mathbf{H}}}
\nc{\bM}{{\mathbf{M}}}
\nc{\bN}{{\mathbf{N}}}
\nc{\bS}{{\mathbf{S}}}
\nc{\bT}{{\mathbf{T}}}
\nc{\bV}{{\mathbf{V}}}
\nc{\bW}{{\mathbf{W}}}
\nc{\bX}{{\mathbf{X}}}
\nc{\bP}{{\mathbf{P}}}
\nc{\bZ}{{\mathbf{Z}}}

\nc{\sA}{{\mathsf{A}}}
\nc{\sB}{{\mathsf{B}}}
\nc{\sC}{{\mathsf{C}}}
\nc{\sD}{{\mathsf{D}}}
\nc{\sF}{{\mathsf{F}}}
\nc{\sK}{{\mathsf{K}}}
\nc{\sM}{{\mathsf{M}}}
\nc{\sO}{{\mathsf{O}}}
\nc{\sQ}{{\mathsf{Q}}}
\nc{\sP}{{\mathsf{P}}}
\nc{\sV}{{\mathsf{V}}}
\nc{\sW}{{\mathsf{W}}}
\nc{\sZ}{{\mathsf{Z}}}
\nc{\sfp}{{\mathsf{p}}}
\nc{\sr}{{\mathsf{r}}}
\nc{\sfb}{{\mathsf{b}}}
\nc{\sfc}{{\mathsf{c}}}
\nc{\sd}{{\mathsf{d}}}
\nc{\sg}{{\mathsf{g}}}
\nc{\sfl}{{\mathsf{l}}}

\nc{\BK}{{\bar{K}}}

\nc{\tA}{{\widetilde{\mathbf{A}}}}
\nc{\tB}{{\widetilde{\mathcal{B}}}}
\nc{\tg}{{\widetilde{\mathfrak{g}}}}
\nc{\tG}{{\widetilde{G}}}
\nc{\TM}{{\widetilde{\mathbb{M}}}{}}
\nc{\tO}{{\widetilde{\mathsf{O}}}{}}
\nc{\tU}{{\widetilde{\mathfrak{U}}}{}}
\nc{\TZ}{{\tilde{Z}}}
\nc{\tZ}{\widetilde{Z}{}}
\nc{\tx}{{\tilde{x}}}
\nc{\tbv}{{\tilde{\bv}}}
\nc{\tfP}{{\widetilde{\mathfrak{P}}}{}}
\nc{\tz}{{\tilde{\zeta}}}
\nc{\tmu}{{\tilde{\mu}}}

\nc{\td}{\ddot{\underline{d}}{}}
\nc{\tzeta}{\widetilde{\zeta}{}}
\nc{\hd}{{\widehat{\underline{d}}}}
\nc{\hG}{{\widehat{G}}}
\nc{\hBP}{\widehat{\mathbb P}{}}
\nc{\hQ}{{\widehat{Q}}}
\nc{\hsM}{\widehat{\mathsf M}{}}
\nc{\hfM}{\widehat{\mathfrak M}{}}
\nc{\hCP}{\widehat{\mathcal P}{}}
\nc{\hCR}{\widehat{\mathcal R}{}}
\nc{\hCS}{{\widehat{\mathcal S}}}
\nc{\hfZ}{\widehat{\mathfrak Z}{}}

\nc{\urho}{\underline{\rho}}
\nc{\uB}{\underline{B}}
\nc{\uC}{{\underline{\mathbb{C}}}}
\nc{\ui}{\underline{i}}
\nc{\ofP}{{\overline{\mathfrak{P}}}}

\nc{\hrho}{{\hat{\rho}}}

\nc{\unl}{\underline}
\nc{\ol}{\overline}
\nc{\one}{{\mathbf{1}}}
\nc{\two}{{\mathbf{t}}}

\nc{\Tot}{{\mathop{\operatorname{\rm Tot}}}}
\nc{\Hilb}{{\mathop{\operatorname{\rm Hilb}}}}
\nc{\CHom}{{\mathop{\operatorname{{\mathcal{H}}\it om}}}}
\nc{\defi}{{\mathop{\operatorname{\rm def}}}}
\nc{\length}{{\mathop{\operatorname{\rm length}}}}

\nc{\Cliff}{{\mathsf{Cliff}}}
\nc{\Fl}{{\mathsf{Fl}}}
\nc{\Fib}{{\mathsf{Fib}}}
\nc{\Coh}{{\mathsf{Coh}}}
\nc{\FCoh}{{\mathsf{FCoh}}}

\nc{\reg}{{\text{\rm reg}}}

\nc{\cplus}{{\mathbf{C}_+}}
\nc{\cminus}{{\mathbf{C}_-}}
\nc{\cthree}{{\mathbf{C}_*}}
\nc{\Qbar}{{\bar{Q}}}
    
\nc{\bh}{{\bar{h}}}
\nc{\bOmega}{{\overline{\Omega}}}
\nc\tGr{\widetilde{\Gr}}

\nc{\seq}[1]{\stackrel{#1}{\sim}}
\nc\ogu{\overline{G/U}}
\nc\chlam{\check{\lam}}

\nc\St{\operatorname{St}}

\nc\uS{\underline{S}}
\nc\QM{\mathcal{QM}}
\nc\FT{\mathsf{FT}}
\nc{\Ca}{\underline{C_a}}
\nc{\SCa}{\underline{S^mC_a}}
\nc{\Cam}{\underline{C_a^{[m]}}}
\nc{\sCa}{\Ca\x_\CA\underline{S^{m-1}C_a}}
\nc{\CCam}{\Ca\x_\CA\underline{C_a^{[m-1]}}}

\begin{document}

\title{Fourier transform 
on the locus of cyclic spectral curves in the Hitchin base}
\author{Andrei Ionov}
\address{
Boston College, Department of Mathematics, Maloney Hall, Fifth Floor,
Chestnut Hill, MA
02467-3806, United States }
\email{ionov@bc.edu}
\begin{abstract}
We compute the Fourier transform of some of the summands of the push-forward of the constant sheaf under the Hitchin map for $\SL_n$ restricted to the locus of cyclic spectral curves inside the Hitchin base (for $\SL_2$ all spectral curves are cyclic) and give an estimate on the support of the Fourier transforms of the other summands.
\end{abstract}

\maketitle
\sec{}{Introduction}


Ever since being introduced in \cite{H} the Hitchin fibration $$\chi\colon\CM\to\CA$$ of the moduli space of semi-stable Higgs $G$-bundles on a smooth compact algebraic curve $C$ of genus $g>1$ over the Hitchin base played a vital role in different areas of mathematics and mathematical physics. One of the most prominent applications is its use in the proof of the fundamental lemma (\cite{N2}). In this and other applications it is important to study the push-forward of the constant sheaf along $\chi$ or, in other words, to understand how the cohomologies of the fibers of $\chi$ vary along the base.

In this paper we approach the study of $\chi_*\uC$ via its Fourier transform in the case of $G=\SL_n$. Unfortunately, although $\CA$ is isomorphic to an affine plane, it could not be canonically identified with a vector space (for $n>2$) and, therefore, at a first glance it is unclear what Fourier transform even means. For this reasons we restrict ourselves to a locus of cyclic spectral curves inside the Hitchin base $\CA$, which is canonically identified with the vector space $H^0(K^n)$ (for $n=2$ it coincides with the whole $\CA$). We hope, however, that with the constructions of the upcoming paper \cite{LN} this obstacle could be overcome and results could be extended to whole $\CA$, as long as there is a suitable definition of what the Fourier transform on $\CA$ is. 

Thus, we study $\FT\chi_*\uC|_{H^0(K^n)}$. For $m<(n-1)(g-1)$ or $m>(2n+1)(n-1)(g-1)$ we completely describe the Fourier transform of the constituents of $R^m\chi_*\uC|_{H^0(K^n)}$ with full support (\reft{1} 1), 2)) as direct sums of IC-extensions of certain simple local systems on the open parts $\Sec^j(\calC_n)-\Sec^{j-1}(\calC_n)$ of secant varieties of $n$-canonical embedding $\calC_n$ of $C$ with $j\le m$ and $m-j$ even. For general $m$ we prove that the support of $\FT R^m\chi_*\uC|_{H^0(K^n)}$ is contained in $\Sec^m(\calC_n)$ (\reft{1} 3)). Finally, for $n=2$ we compute the Fourier transform of the direct summands of $\chi_*\uC$, which do not have full support, known as endoscopic summands (\reft{2}). Our results are painlessly transferred to the case of $K$ being replaced with a general ample line bundle, with the bound $(n-1)(g-1)$ being replace with the half of the maximal number of points $n$-th power of this bundle separates plus 1.

Some similarities with the setting of \cite{D} are worth mentioning, but the results seem to be not related directly. Finally, it would be interesting to understand if our result has any number theoretic applications.

\ssec{}{Acknowledgements} The author is grateful to R.\,Bezrukavnikov for suggesting the problem and to R.\,Bezrukavnikov and Z.\,Yun for helpful discussions.

\ssec{}{Conventions and notations}
We will assume the base field to be $\BC$ unless otherwise stated. All homology and cohomology groups are also assumed to be taken with $\BC$-coefficients. 

For an algebraic variety $X$ we denote by $D^b_c(X) $ the bounded  derived category of constructible sheaves and by $\Perv(X)$ an abelian subcategory of perverse sheaves. 
All the inverse and direct image functors are assumed to be derived unless the opposite is specified.

If $Y\subset X$ is a nonempty locally closed subset and $L$ is a local system on $Y$ we put $\IC(Y, L)$ for the IC-extension of $L$ shifted to the degree $-\dim Y$. To simplify the notations for $Y\subset X$ we denote $\IC(Y, \uC)$ simply by $\uC_Y$. We sometimes put simply $\uC$ for $\uC_X\in \Perv(X)$ if it cannot provide ambiguity.


For a vector space (or a vector bundle) $V$ we denote by $\FT=\FT_V\colon D^b_c(V) \to D^b_c(V^\vee)$ the Fourier transform functor. It preserves the subcategory of perverse sheaves. For more detail and all the properties we will use, see for example \cite{KW}. 


Let $p_n\colon \BA^1-0\to\BA^1-0$ be the $n$-th power map. The group of $n$-th roots of unity $\mu_n$ act on $p_{n,*}\uC$. For $a\in\BZ/n$ let $\LL_a$ be the local system, which is the direct summand of $p_{n,*}\uC$ corresponding to $a$ viewed as a character of $\mu_n$.  
Put $j_0:\BA^1-0\to\BA^1$ for the natural inclusion. If $a\ne 0$ we have $j_{0*}\LL_a=j_{0!}\LL_a=j_{0!*}\LL_a$ and $\FT j_{0!*}\LL_a=j_{0!*}\LL_a$. For $a=0$ we have $\LL_0=\uC_{\BA^1-0}$. We have $j_{0!*}\uC_{\BA^1-0}=\uC_{\BA^1}$ and $\FT\uC_{\BA^1}=\uC_0$. Moreover, $j_{0*}\uC_{\BA^1-0}$ and $j_{0!}\uC_{\BA^1-0}$ are respectively extensions of $\uC_0$ to $\uC_{\BA^1}$ and vice versa and we have $\FT j_{0*}\uC_{\BA^1-0}=j_{0!}\uC_{\BA^1-0}$.

All the quotients by the group actions are understood in a naive sense (not as stacks).

\sec{prel}{Preliminaries}

\ssec{}{Hitchin fibration}

Put $G=SL_n$ and let $\fg=\fsl_n$ be its Lie algebra. We fix a maximal torus $T\subset G$ and a Cartan subalgebra $\ft\subset\fg$. 
We put $W\simeq S_n$ for the Weyl group. The quotient $\ft/W$ is identified with the space of of degree $n$ polynomial of the form $x^n+a_1x^{n-2}+\ldots+a_{n-1}$. In particular, there is is a non-canonical isomorphism $\ft/W\simeq\BA^n$.


 Fix throughout $C$, a smooth connected projective curve of genus $g > 1$, and put $K$ for the canonical line bundle over $C$. 
 We denote by $\CM$ the moduli space of semi-stable Higgs $G$-bundles over $C$, i.e. pairs $(E,\phi)$, where $E$ is a vector bundle of rank $n$ with the trivial determinant and $\phi\in H^0(C, \End(E)\ten K)$ is a Higgs field.

The {\itshape  Hitchin base} $\CA$ is the space of sections of $(K\ten\ft)/W$, which we also identify with the space of degree $n$ polynomial functions from $H^0(K)$ to $H^0(K^n)$ with the coefficient of $x^n$ being $1$ and of $x^{n-1}$ being $0$. In particular, there is a non-canonical isomorphism $$\CA\simeq\bigoplus_{i=2}^n H^0(K^i),$$ where we consider on the right hand side we identify the vector space with the corresponding affine space. 
There is  a "characteristic polynomial of $\phi$" map $$\chi\colon\CM\to\CA,$$ which is called the {\itshape Hitchin map}. For more details, see for example \cite{H}.

For $a\in \CA$ there is a spectral curve $C_a$ defined as a subvariety of the total space of $K$ given by zeroes of the polynomial corresponding to $a$. The natural map $\pi_a\colon C_a\to C$ is a ramified cover of degree $n$. 
We define $\CA_{reg}\subset\CA$ to be the open dense subvariety of $a\in\CA$ for which $C_a$ is smooth and connected. 
 It is known (\cite{BNR}) that for $a\in \CA_{reg}$ we have $\chi^{-1}(a)\cong \Prym_a$, where the Prym variety $\Prym_a\subset \Pic(C_a)$ is the kernel of the norm map $\Pic(C_a)\to \Pic(C)$ induced by $\pi_a$. 
The isomorphism is established by the direct image map $\pi_{a*}$ together with a twist to adjust the determinant. 

For $n=2$ the spectral curve $C_a$ is reduced if $a\ne 0$ and it is reducible if and only if $a$ is a square of a section of a line bundle over $C$ (\cite[Lemma 11.3]{N1}).

\ssec{}{Locus of cyclic covers} Note that, despite the non canonicity of the direct sum presentation, the summand $H^0(K^n)\subset\CA$ is a canonically defined closed subvariety. We can describe it by introducing weights of the variables by making weight of $K$ to be $1$, then quotient the ring of function $\BC[\CA]$ by the ideal generated by all elements of weights strictly between $0$ and $n$. Note that for $n=2$ we have $H^0(K^2)=\CA$. We denote $H^0(K^n)_{reg}=H^0(K^n)\cap\CA_{reg}$ inside $\CA$.

 For $a\in H^0(K^n)$ the spectral curve $C_a$ is given by $x^n+a$ and $\pi_a$ is the cyclic cover. There is an action of a group of $n$-th roots of unity by automorphisms of the cover. For an $n$-th root of unity $\zeta\in\mu_n$ let $\iot_\zeta\colon C_a\to C_a$ be the corresponding automorphism, so that $\pi_a\iot_\zeta=\pi_a$. On cohomology we get a direct sum decomposition with respect to the characters of $\mu_n$-action
 $$H^1(C_a)=\bigoplus_{i\in\BZ/n} H^1(C_a)_i,$$ where $H^1(C_a)_0=H^1(C)$ and $\zeta\in\mu_n$ acts by $\zeta^i$ on $H^1(C_a)_i$.
 
 Let $\Ca$ be the universal spectral curve over $H^0(K^n)$ which comes together with the map $\pi\colon\Ca\to C\x H^0(K^n)$. We denote by $\tilchi\colon\Ca\to H^0(K^n)$ the composition $\proj\circ \pi$ of $\pi$ with the projection on the second factor $\proj\colon C\x H^0(K^n)\to H^0(K^n)$.
 Since $K^n$ is very ample for $n,g\ge2$ the universal family $\Ca$ is smooth. We then have a decomposition with respect to the characters of $\mu_n$-action $$\pi_*\uC_{\Ca}=\bigoplus_{i\in\BZ/n} \bM_i$$ into a direct sum of perverse sheaves. We have  $\bM_0=\uC_{C\x H^0(K^n)}$ and for $i\ne0$ the sheaf $\bM_i$ is a simple perverse sheaf which is an extension by $0$ of a rank $1$ local system away from ramification locus of $\pi$.

\ssec{}{On symmetric powers of curves}

For an algebraic variety $X$ we put $S^mX:= X^m/S_m$ for its $m$-th symmetric power. If $X$ is a smooth curve $S^mX$ is also smooth and is of dimension $m$ for all $m$. 

On $C^{m}\x H^0(K^n)$ consider a sheaf $\bM_i^{\boxtimes m}$ and shift it to the perverse degree. It is a simple perverse sheaf, which has a natural $S_m$-equivariant structure. Note that $S_m$-stabilizers of points in $C^m$ act trivially on  the stalks of $\bM^{\boxtimes m}_i$. It then descends to a simple perverse sheaf $\bM^m_i$ on the quotient $S^mC\x H^0(K^n)$ (see \cite{BL}). It is also Verdier self-dual and is an extension by $0$ from outside the collections of points intersecting ramification locus of $\pi$.

Consider the projection $\proj^m\colon S^mC\x H^0(K^n)\to H^0(K^n)$ on the second fiber. Let $\bM_{i,a}^m$ be a restriction of $\bM^m_i$ to $(\proj^m)^{-1}(a)$ (shifted to the correct degree). For the rest of the subsection we assume that $a\in H^0(K^n)_{reg}$. Then by the same argument the sheaf $\bM^m_{i,a}$ is also a simple perverse. 

The formula of \cite{M} for the cohomology of the symmetric power of a smooth connected algebraic curve $C'$ provides a canonical isomorphism $$H^*(S^mC')\cong \bigoplus_{i+j\le n, i,j\ge0} \left(\bigwedge\nolimits^i H^1(C')\right)[-i-2j],$$ where the index $j$ correspond to the multiplicity of the generator of $H^2(C')$ (hence the shift $-2j$). In particular, the natural map
$$H^*(S^mC)\to H^*(C^m)=(H^*(C))^{\otimes m}$$
is precisely the inclusion of graded $S_n$-invariants. 
The following statement could be considered its refinement with respect to the cyclic covers. 

\lem{coh}
We have $$H^0(S^mC,\bM^m_{i,a})\cong \bigwedge\nolimits^m H^1(C_a)_i$$ and all other cohomology groups of $\bM_{m,a}$ vanish.

\elem

\prf

Let $\Gamma=S_m \wr \mu_n$ be the wreath product group. It acts on $C_a^m$ with the quotient $S^mC$. Denote by $q$ the quotient map. Consider $\uC_{C_a^m}$ as a $\Gam$-equivariant perverse sheaf. Then $q_*\uC_{C_a^m}$ is a also $\Gamma$-equivariant perverse sheaf. Since $\Gamma$-acts trivially on  $S^mC$ the sheaf $q_*\uC_{C_a^m}$ decomposes into a direct sum of components corresponding to irreducible $\Gam$-representations. Note that the generic stalk of $q_*\uC_{C_a^m}$ is a regular representation. We see that $\bM^m_{i,a}$ is precisely a summand corresponding to the one dimensional representation given by character $i$ on each copy of $\mu_n$ and trivial on $S_n$. Computing the cohomology of $q_*\uC_{C_a^m}$ commutes with the direct sum decomposition. Hence, $H^*(S^mC,\bM^m_{i,a})$ is the component of $H^*(q_*\uC_{C_a^m})=H^*(C_a^m)=(H^*(C_a))^{\otimes n}$ corresponding to the one dimensional representation above. But this component is precisely $\bigwedge\nolimits^m H^1(C_a)_i$.




\epr

\ssec{sec}{Secant varieties of $n$-canonical embedding}

For a point $c\in C$ we identify a line, which is a preimage in $H^0(K^n)^\vee$ of an image of $c$ in $\BP(H^0(K^n)^\vee)$ under the $n$-canonical embedding, with the dual to the fiber $K^n_c$ of $K^n$ at $c$, which is embedded by a map $(K^n_c)^\vee\xhookrightarrow{}H^0(K^n)^\vee$ dual to the evaluation at a point map $H^0(K^n)\to  K^n_{c}.$ We will need the following.
\lem{surj}
For a finite collection $c_1,\ldots,c_N$ of pairwise different points on $C$ the evaluation map $e\colon H^0(K^n)\to \bigoplus_{i} K^n_{c_i}$ is surjective if $N<2(n-1)(g-1)$. 
\elem
\prf
Consider a short exact sequence of sheaves
$$0\to K^n(-c_1-\ldots-c_N)\to K^n\to \CO_{c_1}\oplus\ldots\oplus\CO_{c_N}\to0.$$
It yields a long exact sequence of cohomology groups.
$$0\to H^0(K^n(-c_1-\ldots-c_N))\to H^0(K^n)\to$$ $$\to \bigoplus_i K^n_{c_i}\to H^1(K^n(-c_1-\ldots-c_N))\to0,$$
where the map in the middle is $e$. Moreover, $H^1(K^n(-c_1-\ldots-c_N))\cong H^0(K^{-n+1}(c_1+\ldots +c_N))^\vee$. The latter cohomology group is $0$ under the assumptions of the Lemma for degree reasons. 
\epr
\rem{}
The bound is exact, since, if $N=2(n-1)(g-1)$ and the divisor $c_1+\ldots+c_N$ is in the linear system of $K^{n-1}$ the cokernel of $e$ is one-dimensional.

\erem

Let $\calC_n$ be the cone over the image of the $n$-canonical embedding $C\xhookrightarrow{}\BP(H^0(K^n)^\vee)$ and let $\Sec^m(\calC_n)$ be the $m$-th secant variety of $\calC_n$, i.e. the union of all planes spanned by $m$ points on $\calC_n$. In other words a collection of $m$ points on $c_1,\ldots,c_m$ and the elements of the respective fibers $(K^n_c)^\vee$ define a point of $\Sec^m(\calC_n)$. It follows from the Lemma that such presentation is unique for any point in $\Sec^m(\calC_n)-\Sec^{m-1}(\calC_n)$ for $m< (n-1)(g-1)$. 

We also need to study the addition map $$\Sig_{i,j}\colon(\Sec^i(\calC_n)-\Sec^{i-1}(\calC_n))\times(\Sec^j(\calC_n)-\Sec^{j-1}(\calC_n))\to\Sec^{i+j}(\calC_n).$$ 
Assume $i+j<(n-1)(g-1)$. For $k\ge 0$ consider the fiber of $\Sig_{i,j}$ over a point in $\Sec^{i+j-k}(\calC_n)-\Sec^{i+j-k-1}(\calC_n)$ given by a collection $c_1,\ldots, c_{i+j-k}$ of different points on $C$ and nonzero elements $\xi_{c_u}$ of the fibers $(K^n_{c_u})^\vee$. A point of this fiber could be described as follows.

 Pick $0\le l\le\lfloor \frac{k}{2}\rfloor$ and choose $l$ points $\tilc_1,\ldots,\tilc_l$ on $C$ different from $c_1,\ldots, c_{i+j-k}$ and nonzero elements of the fiber $\xi_{\tilc_u}$ at each of them. Separate $c_1,\ldots, c_{i+j-k}$ into three disjoint subsets of sizes $k-2l, i-k+l$ and $j-k+l$ for each point $c_u$ in the first group chose a nonzero element of the fiber at this point $\xi'_{c_u}\ne\xi_{c_u}$. 
Now put  $\tilc_1,\ldots,\tilc_l$ with respective $\xi_{\tilc_u}$ in the fiber, points in the first group with respective $\xi'_{c_u}$ in the fiber and points in the second group with respective $\xi_{c_u}$ in the fiber together to get an element of $\Sec^i(\calC_n)-\Sec^{i-1}(\calC_n)$ and $\tilc_1,\ldots,\tilc_l$ with respective $-\xi_{\tilc_u}$ in the fiber, points in the first group with respective $\xi_{c_u}-\xi'_{c_u}$ in the fiber and points in the third group with respective $\xi_{c_u}$ in the fiber together to get an element of $\Sec^j(\calC_n)-\Sec^{j-1}(\calC_n)$. The resulting point of $\Sec^i(\calC_n)-\Sec^{i-1}(\calC_n))\times(\Sec^j(\calC_n)-\Sec^{j-1}(\calC_n))$ goes to our starting point of $\Sec^{i+j-k}(\calC_n)-\Sec^{i+j-k-1}(\calC_n)$ and each point in the fiber is obtained in this way. 

Choices of $l$ and the ways to separate $c_1,\ldots, c_{i+j-k}$ into three subsets give different irreducible components of this fiber. For a given $l$ there are $\binom{i+j-k}{k-2l,i-k+l,j-k+l}$ of them and each has dimension $k$.




\sec{}{Sheaves on Hitchin base} 

\ssec{full}{Sheaves with full support}

Restricted to $H^0(K^n)_{reg}$ each $\proj_* \bM_i$ is a local system, which is in a single degree for $i\ne0$ and is just constant sheaf for $i=0$. We denote this local systems by $\CL_i$ for $i\ne0$. 
We put $\CL=\bigoplus_{i\ne0}\CL_i$. 
Note also that $\CL_i\simeq\CL_{n-i}^\vee$ with the isomorphism provided by a choice of generator in $H^2(C_a)$.

Over $H^0(K^n)_{reg}$ the isomorphism $\chi^{-1}(a)\cong \Prym_a$ implies an identification   $$R^m\chi_* \uC|_{H^0(K^n)_{reg}}\simeq \bigwedge\nolimits^m\CL=\bigoplus_{m_1+\ldots+m_{n-1}=m}\bigwedge\nolimits^{m_1}\CL_1\ten\ldots\ten\bigwedge\nolimits^{m_{n-1}}\CL_{n-1}.$$  




Recall the relative Hard Lefschetz Theorem for $\chi$. 
Under the decomposition into sum of $\CL_i$ it  could be stated as follows. For each $i$ there is an element in $\CL_i\otimes \CL_{n-i}$, which is skew symmetric if $i=n/2$. It induces maps 
$$\bigwedge\nolimits^{m_i}\CL_i\ten\bigwedge\nolimits^{m_{n-i}}\CL_{n-i}\to \bigwedge\nolimits^{m_i+1}\CL_i\ten\bigwedge\nolimits^{m_{n-i}+1}\CL_{n-i}$$ for $i\ne n/2$ and $$\bigwedge\nolimits^{m_{n/2}}\CL_{n/2}\to \bigwedge\nolimits^{m_{n/2}+2}\CL_{n/2},$$ which are injections for $m_{n-i}+m_{n-i}<\rk \CL_i$ or $m_{n/2}<\rk\CL_{n/2}/2$ and surjections otherwise and in either case they split off as a direct summands. Moreover, composing these maps we get an isomorphism 
$$\bigwedge\nolimits^{m_i}\CL_i\ten\bigwedge\nolimits^{m_{n-i}}\CL_{n-i}\iso \bigwedge\nolimits^{\rk \CL_i -m_{n-i}}\CL_i\ten\bigwedge\nolimits^{\rk \CL_i-m_i}\CL_{n-i}$$ and $$\bigwedge\nolimits^{m_{n/2}}\CL_{n/2}\iso \bigwedge\nolimits^{\rk\CL_{n/2}-m_{n/2}}\CL_{n/2}.$$





\th{1}
1) For $m<(n-1)(g-1)$ the Fourier transform $\FT_{H^0(K^n)_{reg}}\IC(H^0(K^n)_{reg},\bigwedge^m\CL_i)$ is an IC-extensions of a rank $1$ local systems on $\Sec^m(\calC_n)-\Sec^{m-1}(\calC_n)$ for $i\ne n/2$ and is a direct sum of IC-extensions of one rank $1$ local systems on each $\Sec^j(\calC_n)-\Sec^{j-1}(\calC_n)$ with $j\equiv m \pmod{2}$ for $i= n/2$. 

2) More generally, for $m=m_1+\ldots+m_{n-1}<(n-1)(g-1)$ the Fourier transform of the direct summand of $\IC(H^0(K^n)_{reg},\bigwedge\nolimits^{m_1}\CL_1\ten\ldots\ten\bigwedge\nolimits^{m_{n-1}}\CL_{n-1})$ not coming from smaller exterior powers via the relative Hard Lefschetz Theorem is the IC-extension of an irreducible local system on $\Sec^m(\calC_n)-\Sec^{m-1}(\calC_n)$ of rank $\binom{m}{m_1,\ldots,m_{n-1}}$.

3) For any $m$ the support of $\FT_{H^0(K^n)}\IC(H^0(K^n)_{reg},\bigwedge^m\CL)$ is contained inside $\Sec^m(\calC_n)$.

\eth

\rem{}
1) We see, in particular, that for $m<(n-1)(g-1)$ the direct summand of $\bigwedge\nolimits^{m_1}\CL_1\ten\ldots\ten\bigwedge\nolimits^{m_{n-1}}\CL_{n-1}$ not coming from smaller exterior powers via the relative Hard Lefschetz Theorem is in fact an irreducible local system.

2) Using relative Hard Lefschetz isomorphism we extend results of part 1) and 2) to $m>\rk\CL-(n-1)(g-1)$ as well, and for part 3) if $m>\rk\CL/2$ the support is further contained in $\Sec^{\rk\CL-m}(\calC_n)$

\erem

\prf
1) Note that $\bM^m_i$ is simple self-dual of geometric origin and, thus, the decomposition theorem is applicable to it. It follows from \refl{coh} that $\IC(H^0(K^n)_{reg},\bigwedge^m\CL_i)$ is a direct summand in $\proj^m_*\bM^m_i$. 

Since the Fourier transform commute with proper push-forwards we know that $\FT_{H^0(K^n)}\IC(H^0(K^n)_{reg},\bigwedge^m\CL_i)$ is a direct summand of $\FT_{H^0(K^n)} \proj^m_*\bM^m_i=\proj^{\vee m}_*\FT_{S^mC\x H^0(K^n)}\bM^m_a$, where $\proj^{\vee m}\colon S^mC\x H^0(K^n)^\vee\to H^0(K^n)^\vee$ is the projection. Let us compute the last Fourier transform pointwise over $S^mC$.

Fix a partition of $m$, which we represent as a collection of nonnegative integer numbers $(\nu_1, \ldots, \nu_k)$ with $\sum i\nu_i=m$. Consider a point $\uc\in S^mC$, which is a collection of points $c^{(1)}_1, \ldots, c^{(1)}_{\nu_1}\in C$ repeated one time each, points $c^{(2)}_1, \ldots, c^{(2)}_{\nu_2}\in C$ repeated two times each, $\ldots$ , $c^{(k)}_1, \ldots, c^{(k)}_{\nu_k}\in C$ repeated $k$ times each, with all points being pairwise different.

Let $e_\uc\colon H^0(K^n)\to \bigoplus_{\alp,\bet} K^n_{c^{(\bet)}_\alp}$ be the evaluation map, where $K_{c}$ is the fiber of $K$ at the point $c\in C$. By \refl{surj}, the map $e_\uc$ is a surjection. 
By construction of $\bM^m_a$ we have $$\bM^m_i|_{\uc\x H^0(K^n)}=e_\uc^*((j_{0!}\LL_{i})^{\boxtimes \nu_1}\boxtimes (j_{0!}\LL_{2i})^{\boxtimes \nu_2}\boxtimes\ldots\boxtimes (j_{0!}\LL_{ki})^{\boxtimes \nu_k}).$$
Consequently we have $$\FT_{H^0(K^n)}\bM^m_i|_{\uc\x H^0(K^n)}=e^\vee_{\uc!}((j_{0*}\LL_{i})^{\boxtimes \nu_1}\boxtimes (j_{0*}\LL_{2i})^{\boxtimes \nu_2}\boxtimes\ldots\boxtimes (j_{0*}\LL_{ki})^{\boxtimes \nu_k}),$$ where $e^\vee_{\uc}\colon \bigoplus_{\alp,\bet} (K^n_{c^{(\bet)}_\alp})^\vee\to H^0(K^n)^\vee$ is the dual coevaluation map. Note that $\LL_{di}$ is trivial if and only if $d$ is divisible by $n/\mathrm{lcd}(n,i)$.

For $j\le m$ pick a point in $\Sec^j(\calC_n)-\Sec^{j-1}(\calC_n)$, which is described as a collection of distinct points $c_1,\ldots, c_j$ of $C$ and nonzero element of the fibers $(K^n_{c_i})^\vee$, and intersect its preimage under $\proj^{\vee m}$ with $\supp \FT_{S^mC\x H^0(K^n)}\bM^m_i$. The point of the intersection could be described by assigning some multiplicities to the points $c_1,\ldots, c_i$ and adding $l$ additional points with multiplicities divisible by $n/\mathrm{lcd}(n,i)$, 
so that the sum of all the multiplicities is $m$  (in particular, $l\le \lfloor\frac{m-j}{n/\mathrm{lcd}(n,i)}\rfloor$) and then taking the corresponding point in $S^mC$. The $H^0(K^n)^\vee$ direction is then defined uniquely as the image of the coevaluation map for this collection of points with fibers of $(K^n)^\vee$ at $c_1,\ldots, c_i$ being the same as the ones we started with and at the new points being $0$.

Note that at each such point the cohomologies of the stalk of $\supp \FT_{S^mC\x H^0(K^n)}\bM^m_i$ are concentrated in the degree range from $-m-j-l$ to $-m-j$, because stalks of $j_{0*}\LL_{di}$ for $di\ne 0$ are in degree $-1$ and stalks of $j_{0*}\LL_{0}$ are in degrees $0$ and $-1$ at $0$ and in degree $-1$ outside $0$. The total dimension of this fiber is $\lfloor\frac{m-j}{n/\mathrm{lcd}(n,i)}\rfloor$ (maximal number of choices of the new points). We see that the cohomologies of the stalks of $\FT_{H^0(K^n)} \proj^m_*\bM_m=\proj^{\vee m}_*\FT_{S^mC\x H^0(K^n)}\bM^m_i$ at points of $\Sec^j(\calC_n)-\Sec^{j-1}(\calC_n)$ are concentrated in degrees below $-m-j+2\lfloor\frac{m-j}{n/\mathrm{lcd}(n,i)}\rfloor\le -2j$, which are nonpositive perverse degrees. The equality is achieved if and only if $m=j$ or $i=n/2$ and $m-j$ is even. In this case all the multiplicities are assigned uniquely: 1 to $c_1,\ldots,c_j$ and $2$ to the new points if needed. Moreover, the top dimension component of the fiber is unique as it contains an open dense subset isomorphic to an open dense subset of $\frac{m-j}{2}$ pairwise different points and the local system restricted to this open subset has rank $1$, so that the cohomology in the perverse degree are at most $1$-dimensional.

We conclude that $\FT_{H^0(K^n)} \proj^m_*\bM_m$ is an irreducible perverse sheaf if $i\ne n/2$ and is a sum of at most $\lfloor\frac{m}{2}\rfloor+1$ simple perverse sheaves if $i=n/2$. Knowing that $\IC(H^0(K^n)_{reg},\bigwedge^m\CL_i)$ is a direct summand in $\proj^m_*\bM^m_i$ and it has at least $\lfloor\frac{m}{2}\rfloor+1$ direct summands if $i=n/2$, we deduce that $\proj^m_*\bM^m_i\cong \IC(H^0(K^n)_{reg},\bigwedge^m\CL_i)$. This proves the theorem.



2) The Fourier transform sends tensor product to the !-convolution. Therefore, $\FT_{H^0(K^n)}\IC(H^0(K^n)_{reg},\bigwedge\nolimits^{m_1}\CL_1\ten\ldots\ten\bigwedge\nolimits^{m_{n-1}}\CL_{n-1})$ is a constituent of $\FT_{H^0(K^n)}\IC(H^0(K^n)_{reg},\bigwedge^{m_1}\CL_1)\star\ldots\star\FT_{H^0(K^n)}\IC(H^0(K^n)_{reg},\bigwedge^{m_{n-1}}\CL_{n-1})$. We will now compute the latter sheaf using that we know the factors from part 1). Assume that for $i=n/2$ we take only the direct summand with the maximal support. 

It follows from the description of \refss{sec} that the convolution map is semismall. We, thus, need to just track the top degree cohomology of each of the irreducible components. Note that $\FT_{H^0(K^n)}\IC(H^0(K^n)_{reg},\bigwedge^{m_i}\CL_i)$ is $\LL_i$ restricted to any line contained in $\Sec^{m_i}(\calC_n)-\Sec^{m_i-1}(\calC_n)$ and passing through zero. We see that over $\Sec^{m-k}(\calC_n)-\Sec^{m-k-1}(\calC_n)$ only the components corresponding to adding $k/2$ new points have nonzero cohomology in the top degree (in particular, $k$ is even) and only if the new points are added to the opposite summands $i$ and $n-i$, as only in this case the local system restricted to the component has top degree cohomology. In the latter case the top degree cohomology are one-dimensional. The number of such components match with the ranks of the summands coming from relative Hard Lefschetz. We conclude that the Fourier transform of the summand not coming from smaller exterior powers via the relative Hard Lefschetz Theorem is supported on $\Sec^{m}(\calC_n)-\Sec^{m-1}(\calC_n)$.

Finally, there are exactly $\binom{m}{m_1,\ldots,m_{n-1}}$ points over the generic point of $\Sec^m(\calC_n)$. It remains to show that the corresponding local system is irreducible. Pick a vector in a stalk of the local system over a point given by a collection  of different points 
of $C$ and nonzero elements 
of the respective fibers of $(K^n)^\vee$. 
By going around the loop that interchanges two points of $C$ and taking linear combination we can achieve that this vector has nonzero entries in all positions in the standard basis corresponding to the points in the fiber of the addition map.  Now since the monodromy action around the loop moving along the fiber of $(K^n)^\vee$ at a given point around $0$ acts differently on different summands we can span the whole stalk, thus, proving that the local system is irreducible.

3) This essentially follows from the arguments of the first two parts. By the argument of part 1) $\FT_{H^0(K^n)} \proj^m_*\bM^m_i$ is still supported inside $\Sec^m(\calC_n)$ regardless of $m$ and the convolution of sheaves on secant varieties $\Sec^{m_i}(\calC_n)$ is supported inside the secant variety $\Sec^m(\calC_n)$ for the part 2) of the argument. 

\epr

\ssec{red}{Constant sheaf on the endoscopic locus for $n=2$}

In this section we restrict to the case of $n=2$. For a nontrivial line bundle $\CT\in\Pic(C)[2]$ with trivial square $\CT^{\otimes 2}\simeq\CO$ consider a square map $s_\CT\colon H^0(K\otimes \CT)\to H^0(K^2)$. Put $\BS_\CT\subset H^0(K^2)$ for its image. Note that, from the Riemann-Roch theorem it follows immediately that 
$\dim H^0(K\otimes \CT)=g-1$.

By the description of  \cite[Section 11]{N1} the constituents of $\chi_*\uC$ without full support (and not just supported at $0$) are supported on $\BS_\CT$. Furthermore, this constituents are constant along these subvarieties, as they relate to the Hitchin fibrations for the endoscopic group, which in this case is a torus, thus these fibrations are trivial. In this section we compute $\FT_\CA\uC_{\BS_\CT}$.



A point $\xi\in\BP(H^0(K^2)^\vee)$ defines a quadric $Q_{\CT,\xi}\subset\BP(H^0(K\otimes \CT))$ by the equation $\xi(s_\CT(x))=0$. Let $Z_{\CT,i}\subset \BP(H^0(K^2)^\vee)$ be a locally closed subvariety consisting of points $\xi$ such that the quadric $Q_{\CT,\xi}$ has rank $i$ and let $\tilZ_{\CT,i}\subset H^0(K^2)^\vee$ be the cone over $Z_{\CT,i}$ with $0$ removed. Collections $\{Z_{\CT,i}\}$ and $\{\tilZ_{\CT,i}\}\cup \{0\}$ provide the stratifications of $\BP(H^0(K^2)^\vee)$ and $H^0(K^2)^\vee$ respectively.

\th{2}

The Fourier transform $\FT_\CA\uC_{\BS_\CT}$ is an irreducible IC-sheaf, which restricts to $0$ on $\tilZ_{\CT,i}$ if $i$ is odd and is a rank $1$ local system in with monodromy group $\{\pm 1\}$ restricted to $\tilZ_{\CT,i}$ if $i$ is even. 

\eth{}


\prf
The sheaf $\uC_{\BS_\CT}$ is an irreducible IC-sheaf and then so is $\FT_\CA\uC_{\BS_\CT}$.

By \cite{B} it suffices to compute the Radon transform of the constant sheaf on the image of $\BS_\CT$ in $\BP(H^0(K^2))$, whose stalk at $\xi\in \BP(H^0(K^2)^\vee)$ is then the reduced cohomology group $\tilH^*(\BP(H^0(K\ten\CT))-Q_{\CT,\xi})\simeq\Ker(H_*(Q_{\CT,\xi})\to H_*(\BP(H^0(K\ten\CT))))^\vee$. This group is nonzero only if rank of $Q_{\CT,\xi}$ is even and only in one degree equal to the rank, where it is one dimensional with the action of monodromy by $\pm 1$.

\epr

\end{document}